\documentclass[preprint,12pt]{elsarticle}
\usepackage[numbers]{natbib}
\usepackage{rotating}
\usepackage{times}
\usepackage{graphicx,multirow,lineno}
\usepackage{amsmath,amssymb,bbm,amsthm}
\usepackage{pifont}

\def\beq{\begin{equation}}
\def\eeq{\end{equation}}
\def\bd{\begin{description}}
\def\ed{\end{description}}
\newtheorem{theorem}{Theorem}
\newtheorem{lemma}{Lemma}
\def\G1{\hbox{$\displaystyle{\mbox{\ding{172}}}$}}
\journal{Communications in Nonlinear Science and Numerical Simulation}
 
\begin{document}


 \begin{frontmatter}
\title{On strong homogeneity of a class of global optimization algorithms working with infinite and infinitesimal scales}


\author[addDIMES,addNNGU]{Yaroslav~D.~Sergeyev\corref{cor1}}
\ead{yaro@dimes.unical.it} \cortext[cor1]{Corresponding author.}
\author[addDIMES,addNNGU]{Dmitri~E.~Kvasov}
\ead{kvadim@dimes.unical.it}
\author[addDIMES,addNNGU]{Marat~S.~Mukhametzhanov}
\ead{m.mukhametzhanov@dimes.unical.it}
\address[addDIMES]{Dipartimento di Ingegneria Informatica, Modellistica, Elettronica e Sistemistica,\\
Universit\`{a} della Calabria, Rende (CS), Italy\\[2pt]}
\address[addNNGU]{Department of Software and Supercomputing Technologies, \\
Lobachevsky State University, Nizhni Novgorod, Russia}


\small{
\begin{abstract}
The necessity to find the global optimum of multiextremal functions
arises in many applied problems where finding local solutions is
insufficient. One of the desirable properties of global optimization
methods is \emph{strong homogeneity}  meaning that a method produces
the same sequences of points where the objective function is
evaluated independently both of multiplication of the function by a
scaling constant and of adding a shifting constant.   In this paper,
several aspects of global optimization   using  strongly homogeneous
methods are considered.  First, it is shown that even if a method
possesses this property theoretically, numerically very small and
large scaling constants can lead to ill-conditioning of the scaled
problem. Second, a new class of global optimization  problems where
the objective function can have not only finite but also infinite or
infinitesimal Lipschitz constants is introduced. Third, the strong
homogeneity of several Lipschitz global optimization algorithms is
studied in the framework of the Infinity Computing paradigm allowing
one to work \emph{numerically} with a variety of infinities and
infinitesimals. Fourth, it is proved that a class of efficient
univariate methods enjoys this property for finite, infinite and
infinitesimal scaling and shifting constants. Finally, it is shown
that in certain cases the usage of numerical infinities and
infinitesimals can avoid ill-conditioning produced by scaling.
Numerical experiments illustrating theoretical results are
described.
\end{abstract}
}

\begin{keyword} Lipschitz global optimization \sep
strongly homogeneous methods \sep numerical infinities and
infinitesimals \sep ill-conditioned problems
\end{keyword}
\end{frontmatter}
\newpage
\section{Introduction}

 In many applied problems it is required to find the global optimum
(minimization problems are considered here, i.e., we talk about the
global minimum) of multiextremal non-differentiable functions. Due
to the presence of multiple local minima and non-differentiability
of the objective function, classical local optimization techniques
cannot be used for solving these problems and global optimization
methods should be developed (see, e.g.,
\cite{Floudas&Pardalos(2009),Jones&Perttunen&Stuckman(1993),
Paulavicius&Zilinskas(2014),
Pinter(2002),Sergeyev&Kvasov(2017),Strongin(1978),
Strongin&Sergeyev(2000),Zilinskas&Zilinskas(2013),Zhigljavsky&Zilinskas(2008),
Zilinskas(2012)}).

 One of the desirable properties of   global
optimization methods (see
\cite{Elsakov&Shiryaev(2006),Strongin(1978),Zilinskas(2012)}) is
their \emph{strong homogeneity} meaning that a method produces the
same sequences of trial points (i.e., points where the objective
function $f(x)$ is evaluated) independently of both shifting  $f(x)$
vertically and its multiplication   by a scaling constant. In other
words, it can be useful to optimize a scaled function
 \beq
 g(x) = g(x; \alpha, \beta)= \alpha f(x) + \beta, \hspace{1cm} \alpha>0, \label{g}
 \eeq
 instead of the original objective function
$f(x)$. The concept of strong homogeneity has been introduced in
\cite{Zilinskas(2012)} where it has been shown that both the
P-algorithm (see \cite{Zilinskas(1985)}) and   the one-step Bayesian
algorithm (see \cite{Mockus(1988)}) are strongly homogeneous. The
case $\alpha=1, \beta \neq 0$ was considered in
\cite{Elsakov&Shiryaev(2006),Strongin(1978)} where a number of
methods enjoying this property and called    \emph{homogeneous} were
studied. It should be mentioned that there exist   global
optimization  methods that are homogeneous or strongly homogeneous
and algorithms (see, for instance,   the  DIRECT algorithm from
\cite{Jones&Perttunen&Stuckman(1993)} and a huge number of its
modifications) that do not possess this property.

All the  methods mentioned above have been developed to work with
Lipschitz global optimization problems that can be met very
frequently in practical applications (see, e.g.,
\cite{Paulavicius&Zilinskas(2014),Pinter(2002),Sergeyev&Kvasov(2017),Strongin(1978),
Strongin&Sergeyev(2000),Zhigljavsky&Zilinskas(2008)}). These methods
belong to the class of ``Divide-the-best'' algorithms introduced in
\cite{Divide_the_Best}. Efficient methods from this class that
iteratively subdivide the search region and estimate local and
global  Lipschitz constants during the search are studied in this
paper, as well. Two  kinds of algorithms are taken into
consideration: geometric and information ones  (see
\cite{Sergeyev&Kvasov(2017),Sergeyev:et:al.(2013),
Strongin(1978),Strongin&Sergeyev(2000)}).  The first class of
algorithms is based on a  geometrical interpretation of the
Lipschitz condition and takes its origins in the method proposed in
\cite{Piyavskij(1972)} that builds a piece-wise linear minorant for
the objective function using the Lipschitz condition. The second
approach uses a stochastic model developed in \cite{Strongin(1978)}
that allows one to calculate probabilities of locating global
minimizers within each of the subintervals of the search region and
is based on the information-statistical algorithm proposed in
\cite{Strongin(1978)} (for other rich ideas in stochastic global
optimization see \cite{Zhigljavsky&Zilinskas(2008),
Zilinskas(2012)}). Both classes of methods use in their work
different strategies to estimate global and local Lipschitz
constants (see, e.g., \cite{Piyavskij(1972),Sergeyev&Kvasov(2017),
Sergeyev:et:al.(2016a),Sergeyev:et:al.(2013), Strongin(1978),
Strongin&Sergeyev(2000)}).

In this paper, it will be shown that several fast univariate methods
using   local tuning techniques  to accelerate the search through a
smart balancing of the global and local information collected during
the search (see recent surveys in
\cite{Sergeyev&Kvasov(2017),Sergeyev:et:al.(2016a)}) enjoy the
property of the strong homogeneity. In particular, it will be proved
that  this  property   is valid for the considered methods not only
for finite values of the constants $\alpha$ and $\beta$ but for
infinite and infinitesimal ones, as well. To prove this result, a
new class of global optimization problems with the objective
function having infinite or infinitesimal Lipschitz constants is
introduced. Numerical computations with functions that can assume
infinite and infinitesimal values are executed  using the Infinity
Computing paradigm allowing one to work \emph{numerically} with a
variety of infinities and infinitesimals on a patented in Europe and
USA new supercomputer called  the Infinity Computer (see, e.g.,
surveys
 \cite{Sergeyev(2008),Sergeyev(2017)}). This computational
methodology has already been successfully applied in optimization
and numerical differentiation
\cite{Cococcioni,DeLeone,DeLeone_2,Num_dif} and in a number of other
theoretical and applied research areas such as, e.g., cellular
automata \cite{DAlotto},  hyperbolic geometry \cite{Margenstern_3},
percolation \cite{Iudin_2}, fractals \cite{Caldarola_1,Biology},
infinite series \cite{Dif_Calculus,Zhigljavsky}, Turing machines
\cite{Sergeyev_Garro}, numerical solution of ordinary differential
equations \cite{ODE_3,ODE_2}, etc. In particular, in the recent
paper \cite{Gaudioso&Giallombardo&Mukhametzhanov(2018)}, numerical
infinities and infinitesimals from
\cite{Sergeyev(2008),Sergeyev(2017)} have been successfully used to
handle ill-conditioning in a multidimensional optimization problem.

The importance to have the possibility to work with infinite and
infinitesimal scaling/shifting constants $\alpha$ and $\beta$ has an
additional value due to the following fact. It can happen that even
if a method possesses the strong homogeneity property theoretically
and the original objective function $f(x)$ is well-conditioned,
numerically very small and/or large finite constants $\alpha$ and
$\beta$ can lead to the ill-conditioning of the global optimization
problem involving $g(x)$ due to overflow and underflow taking place
when $g(x)$ is constructed from $f(x)$. Thus, global minimizers can
change their locations and the values of global minima can change,
as well. As a result, applying methods possessing the strong
homogeneity property to solve these problems will lead to finding
the changed values of minima related to  $g(x)$ and not the desired
global solution of the original function  $f(x)$ we are interested
in. In this paper, it is shown that numerical infinities and
infinitesimals and the Infinity Computing framework can help in this
situation.

The rest of paper is structured as follows. Section~2 states the
problem formally, discusses ill-conditioning induced by scaling, and
briefly describes the  Infinity Computer framework. It is stressed
that the introduction of numerical infinities and infinitesimals
allows us to consider a new class of functions having infinite or
infinitesimal Lipschitz constants. Section~3 presents geometric and
information Lipschitz global optimization algorithms studied in this
paper and shows how an adaptive estimation of global and local
Lipschitz constants can be performed. So far, the fact whether these
methods are strongly homogeneous or not was an open problem even for
finite constants $\alpha$ and $\beta$. Section~4 proves that these
methods enjoy the strong homogeneity property for finite, infinite,
and infinitesimal scaling and shifting constants. Section~5 shows
that in certain cases the usage of numerical infinities and
infinitesimals can avoid ill-conditioning produced by scaling and
illustrates these results numerically. Finally, Section~6 contains a
brief conclusion.

\section{Problem statement, ill-conditioning induced by scaling,
and the Infinity\\ Computer framework}

\subsection{Lipschitz global optimization and strong homogeneity}

 Let us consider the
following univariate global optimization problem where it is
required to find the global minimum $f^*$ and global minimizers
$x^*$ such that
\begin{equation}
f^* = f(x^*) = \min~f(x), \hspace{5mm} x\in D = [a,b] \subset
\mathbb{R}. \label{problem1}
\end{equation}
It is supposed that the objective function $f(x)$ can be
multiextremal and non-differentiable. Moreover, the objective
function $f(x)$ is supposed to be Lipschitz continuous over the
interval $D$, i.e., $f(x)$ satisfies the following condition
\begin{equation}
|f(x_1) - f(x_2)| \leq L |x_1 - x_2|, \hspace{5mm} x_1,x_2 \in D,
\label{lipschitz1}
\end{equation}
where $L$ is the Lipschitz constant, $0< L < \infty$.

A vast literature  is dedicated to the problem (\ref{problem1}),
(\ref{lipschitz1}) and algorithms for its solving (see, e.g.,
\cite{Floudas&Pardalos(2009),Gergel&Grishagin&Gergel(2016),Grishagin_Israfilov_Sergeyev_2018,
Kvasov&Mukhametzhanov(2018),
Paulavicius&Zilinskas(2014),Sergeyev_Grishagin_1994,Sergeyev&Kvasov(2017),Sergeyev:et:al.(2017b),
Sergeyev:et:al.(2013),Strongin&Sergeyev(2000),
Zhigljavsky&Zilinskas(2008)}). In particular, in practice it can be
useful to optimize a scaled function $g(x)$ from (\ref{g}) instead
of the original objective function $f(x)$ (see, e.g.,
\cite{Elsakov&Shiryaev(2006),Strongin(1978),Zilinskas(2012)}). For
this kind of problems, the concept of strong homogeneity for global
optimization algorithms  has been introduced in
\cite{Zilinskas(2012)}: An algorithm is called \textit{strongly
homogeneous} if it generates the same sequences of trials
(evaluations of the objective function) during optimizing the
original objective function $f(x)$ and the scaled function $g(x)$
from (\ref{g}), where $\alpha>0$ and $\beta$ are   constants (notice
that homogeneous methods corresponding to the case $\alpha=1, \beta
\neq 0$ have been considered originally in
\cite{Elsakov&Shiryaev(2006),Strongin(1978)}). Unfortunately, in
practice it is not always possible to obtain correct values of
$g(x)$ for huge and small values of $\alpha>0$ and $\beta$ due to
overflows and underflows present if traditional computers and
numeral systems are used for evaluation of $g(x)$ even if the
original function $f(x)$ is well-conditioned.

\subsection{Ill-conditioning produced by scaling}

\begin{figure}[t!]
\includegraphics[width = 1\linewidth, height = 110mm]{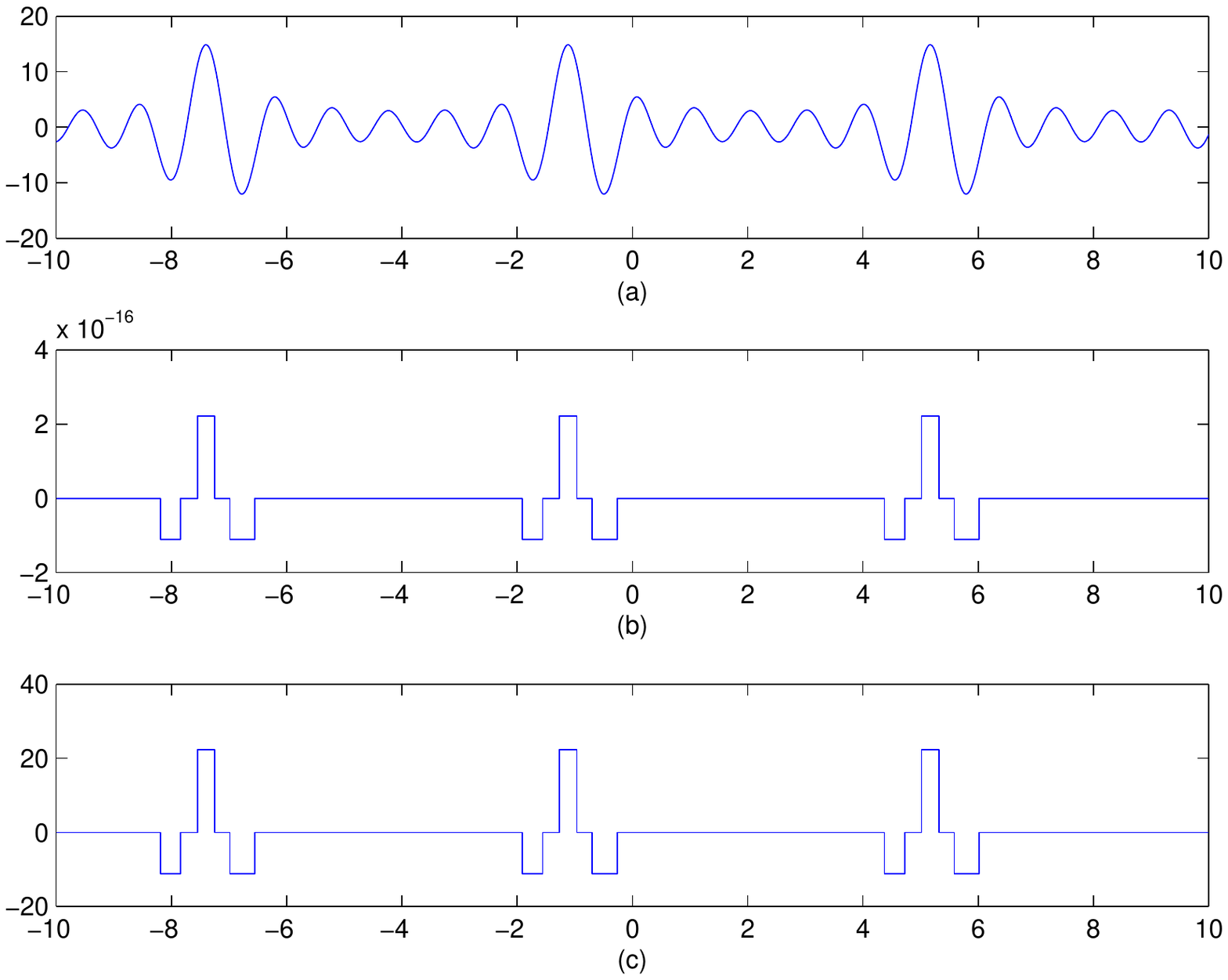}
\caption{Graphs of (a) the test function (\ref{f3}), (b) the scaled
function $g_3(x)$ from (\ref{g3}) in the logarithmic form, (c) the
inverted scaled function $\widehat{f}_3(x) = 10^{17} (g_3(x) - 1)$.
It can be seen that the form of the functions $g_3(x)$ and
$\widehat{f}_3(x)$ are qualitatively different with respect to the
original function $f_3(x)$ due to overflows and underflows.}
\label{fig:f}
\end{figure}

As an illustration, let us consider  the following test problem
from \cite{Hansen&Jaumard(1995)} shown in Fig.~\ref{fig:f}.a:
\begin{equation}
f_3(x) = \sum_{k=1}^5 -k\cdot sin[(k+1)x+k],\hspace{5mm} x \in D =
[-10,10]. \label{f3}
\end{equation}
The function $f_3(x)$ has been chosen from the set of 20 test
functions described in \cite{Hansen&Jaumard(1995)} because it has
the highest number of local minima among these functions and the
following three global minimizers
 \beq
 x_1^* = -0.491, \hspace{5mm}x_2^* = -6.775, \hspace{5mm}x_3^* =
 5.792  \label{minima x}
 \eeq
corresponding to the global minimum
 \beq
 f^* = f(x_1^*) = f(x_2^*) = f(x_3^*) = -12.0312.\label{minima f}
 \eeq
Let us take $\alpha=10^{-17}$ and $\beta=1$ obtaining so the
following function
\begin{equation}
g_3(x) = 10^{-17} f_3(x) + 1. \label{g3}
\end{equation}
It can be seen from Fig.~\ref{fig:f}.a and Fig.~\ref{fig:f}.b that
$f_3(x)$ and $g_3(x)$ are completely different. If we wish to
reestablish $f_3(x)$ from $g_3(x)$, i.e., to compute the inverted
scaled function $\widehat{f}_3(x)= 10^{17} (g_3(x)-1)$, then it will
not coincide with $f_3(x)$. Fig.~\ref{fig:f}.c shows
$\widehat{f}_3(x)$ constructed from $g_3(x)$ using
MATLAB\textsuperscript{\textregistered}  and the piece-wise linear
approximations with the integration step $h = 0.0001$.

 Thus, this scaling leads to an
ill-conditioning. Due to underflows taking place in commonly used
numeral systems (in this case, the type \emph{double} in
MATLAB\textsuperscript{\textregistered}), the function $g_3(x)$
degenerates over many intervals in constant functions and many local
minimizers disappear (see Fig.~\ref{fig:f}.b). In the same time, due
to overflows, several local minimizers become  global minimizers of
the scaled function $g_3(x)$. In particular, using the following two
commands in MATLAB\textsuperscript{\textregistered}
 $$
 [gmin,~imin] = min(y), \hspace{1cm}xmin =
x(imin)
$$
we can calculate an   approximation of the global minimum for
$g_3(x)$. Using the array $y$ containing the values of $g_3(x)$
calculated with the stepsize $h = 0.0001$, i.e.,
 $$y_i = 10^{-17}
f_3(x_i)+1,\hspace{1cm}x_i = -10+h\cdot (i-1),~i\geq 1,
$$
we get   $(xmin,gmin) = (-8.194, 1.0)$  being an approximation of
the global
 minimum
$(x^*,g_3(x^*))$ of $g_3(x)$ that does not coincide with the global
minima (\ref{minima x}), (\ref{minima f}) of the original function
$f_3(x)$. Thus, due to underflows and overflows, the ``wrong''
global minimum of the scaled function $g_3(x)$ has been found.
Analogously, due to the same reasons, the inverted function
$\widehat{f}_3(x) = 10^{17} (g_3(x)-1)$ has also different global
minima with respect to the original function $f_3(x)$ (see
Fig.~\ref{fig:f}.c). Clearly, a similar situation can be observed if
larger values of $\alpha$ and $\beta$ are used (for instance,
$\alpha = 10^{17}$ and $\beta = 10^{35}$).

This example shows that  in   case of very huge or very small finite
values of constants~$\alpha$ and $\beta$, even if it has been proved
theoretically that a method is strongly homogeneous, it does not
make sense to talk about this property since it is not possible to
construct correctly the corresponding scaled functions on the
traditional computers.

\subsection{Infinity Computing briefly}

 The already mentioned Infinity Computing computational
paradigm  (see, e.g., surveys in
\cite{Sergeyev(2008),Sergeyev(2017)}) proposed for working
numerically with infinities and infinitesimals can be used in the
context of strongly homogeneous global optimization methods, as
well. This computational methodology has already been successfully
applied in a number of applications mentioned above. In particular,
it has been successfully used for studying strong homogeneity of the
P-algorithm and the one-step Bayesian algorithm (see
\cite{Zilinskas(2012)}) and
 to handle
ill-conditioning in local optimization (see
\cite{Gaudioso&Giallombardo&Mukhametzhanov(2018)}).

In this paper, it is shown that within the Infinity Computing
paradigm not only finite, but also infinite and infinitesimal values
of $\alpha$ and $\beta$ can be adopted. In particular, the
ill-conditioning present in the global optimization problem
described above in the traditional computational framework can be
avoided in certain cases within the Infinity Computing paradigm.
This is done by using numerical infinite and/or infinitesimal values
of $\alpha$ and $\beta$ instead of huge or very small
scailing/shifting constants. In order to study the strong
homogeneity property with infinite and infinitesimal
scaling/shifting constants, let us introduce the Infinity Computing
paradigm briefly.

 Finite, infinite, and
infinitesimal numbers in this framework are represented using the
positional numeral system with the infinite base $\G1$ (called
\emph{grossone}) introduced as the number of elements of the set of
natural numbers\footnote{It should be emphasized  that the Infinity
Computing approach allows us a full numerical treatment of both
infinite and infinitesimal numbers whereas the non-standard analysis
(see \cite{Robinson}) has a symbolic character and, therefore,
allows symbolic computations only (see a detailed discussion on this
topic in \cite{Sergeyev(2017)}).}. In the \G1-based positional
system a number $C$ expressing the quantity
 \begin{equation}
 C = c_{p_m}\G1^{p_m}+c_{p_{m-1}}\G1^{p_{m-1}}+...+c_{p_1}\G1^{p_1}+
 c_{p_0}\G1^{p_0}+c_{p_{-1}}\G1^{p_{-1}}+...+c_{p_{-k}}\G1^{p_{-k}},
 \end{equation}
 is written in the form
\begin{equation}
C = c_{p_m}\G1^{p_m}...c_{p_1}\G1^{p_1}c_{p_0}\G1^{p_0}c_{p_{-1}}\G1^{p_{-1}}...c_{p_{-k}}\G1^{p_{-k}}.
\label{grossnumber}
\end{equation}
In (\ref{grossnumber}), all numerals $c_i$ are not equal to zero
(they can be positive or negative). They are finite, written in a
traditional numeral system and are called \emph{grossdigits},
whereas all numbers $p_i$, called \emph{grosspowers}, are sorted in
decreasing order with $p_0 = 0$:
\begin{equation}
p_m>...>p_1>p_0>p_{-1}>...>p_{-k}.
\end{equation}
In the \G1-based numeral system, all finite numbers $n_{finite}$
can be represented using only one grosspower $p_0 = 0$ and the
grossdigit $c_0 = n_{finite}$ since $\G1^{0}=1$. The simplest
infinite numbers in this numeral system are expressed by numerals
having at least one finite grosspower greater than zero. Simple
infinitesimals are represented by numerals having only finite
  negative grosspowers. The simplest number from this group is
$\G1^{-1}$ being the inverse element with respect to multiplication
for $\G1$:
 \[
 \frac{1}{\G1} \cdot \G1 = \G1 \cdot \frac{1}{\G1} = 1.
 \]
 It should be mentioned also, that in this framework
numbers having a finite part and infinitesimal ones (i.e., in
(\ref{grossnumber}) it follows $c_j = 0,~j>0,~c_0\neq 0,$ and
$c_i\neq 0$ for at least one $i<0$) are called \emph{finite}, while
the numbers with only one grossdigit $c_0 \neq 0$ and $c_i =
0,~i\neq 0,$ are called \emph{purely finite}. However, hereinafter
all \emph{purely finite} numbers will be called   \emph{finite} just
for simplicity.

\subsection{Functions with infinite/infinitesimal Lipschitz constants}

The introduction of the Infinity Computer paradigm allows us to
consider univariate global optimization problems  with the objective
function $g(x)$ from (\ref{g}) that can assume not only finite
values, but also infinite and infinitesimal ones. It is supposed
that the original function $f(x)$ can assume finite values only and
it satisfies condition (\ref{lipschitz1}) with a finite constant
$L$. However, since in (\ref{g})  the scaling/shifting parameters
$\alpha$ and~$\beta$ can be not only finite, but also infinite and
infinitesimal and, therefore, to work with $g(x)$, the Infinity
Computing framework is required. Thus, the following optimization
problem is introduced
\begin{equation}
 \min~g(x) = \min~(\alpha f(x) + \beta),\hspace{5mm}x\in D = [a,b] \subset
\mathbb{R}, \alpha > 0, \label{problem2}
\end{equation}
where the function $f(x)$ can be multiextremal, non-differentiable,
and Lipschitz continuous with a finite value of the Lipschitz
constant $L$ from (\ref{lipschitz1}). In their turn, the values
$\alpha$ and $\beta$ can be finite, infinite, and infinitesimal
numbers representable in the numeral system (\ref{grossnumber}).

The finiteness of the original Lipschitz constant $L$ from
(\ref{lipschitz1})   is the essence of the Lipschitz condition
allowing people to construct optimization methods for traditional
computers. The scaled objective function $g(x)$ can assume not only
finite, but also infinite and infinitesimal values and, therefore,
in these cases it is not Lipschitzian in the traditional sense.
However, the Infinity Computer paradigm extends the space of
functions that can be treated theoretically and numerically to
functions assuming infinite and infinitesimal values. This fact
allows us to extend the concept of Lipschitz functions to the cases
where the Lipschitz constant can assume infinite/infinitesimal
values.

 Let
us indicate in the rest of the paper by ``~$\widehat{ }$~'' all the
values related to the function $g(x)$ and without ``~$\widehat{
}$~''  the values related to the function $f(x)$.  The following
lemma shows a simple but important property of the Lipschitz
constant for the objective function $g(x)$.
\begin{lemma} The Lipschitz constant $\widehat{L}$ of the function
 $g(x) = \alpha f(x) + \beta$, where $f(x)$ assumes only finite
 values and has the finite Lipschitz constant $L$ over
  the interval $[a,b]$ and $\alpha,~\alpha>0,$ and $\beta$ can
   be finite, infinite, and infinitesimal, is equal to $\alpha L$.
\end{lemma}
\emph{Proof.} The following relation can be obtained  from the
definition of $g(x)$ and the fact that $\alpha>0$
$$|g(x_1)-g(x_2)| = \alpha |f(x_1) - f(x_2)|, \hspace{5mm}
x_1,~x_2 \in [a,b].$$
 Since $L$ is the Lipschitz constant for $f(x)$,
then
$$ \alpha |f(x_1)-f(x_2)| \leq
\alpha L |x_1-x_2| = \widehat{L} |x_1 - x_2|, \hspace{5mm}x_1,~x_2
\in [a,b],$$
 and this inequality proves the lemma. \hfill $\qed$

 Thus, the new Lipschitz condition for the function $g(x)$ from
 (\ref{g}) can be written as
 \beq
  |g(x_1) - g(x_2)| \leq \alpha L |x_1 - x_2| = \widehat{L} |x_1 - x_2|, \hspace{5mm} x_1,x_2 \in D,
\label{lipschitz_g}
 \eeq
where the constant $L$ from (\ref{lipschitz1})  is  finite and the
quantities $\alpha$ and $\widehat{L} $ can assume infinite and
infinitesimal values.

Notice that in the introduced class of functions   infinities and
infinitesimals are expressed in numerals (\ref{grossnumber}), and
Lemma~1 describes the first property of this class. Notice also that
symbol $\infty$ representing a generic infinity cannot be used
together with numerals (\ref{grossnumber}) allowing us to
distinguish a variety of infinite (and infinitesimal) numbers.
Analogously, Roman numerals (I, II, III, V, X, etc.) that do not
allow to express zero and negative numbers are not used in the
positional numeral systems where new symbols (0, 1, 2, 3, 5, etc.)
are used to express numbers.

Some geometric and information global optimization  methods (see
\cite{Pinter(2002),Piyavskij(1972),Sergeyev&Kvasov(2017),
Sergeyev:et:al.(2016a),Sergeyev:et:al.(2013),Strongin(1978),Strongin&Sergeyev(2000)})
used for solving the traditional Lipschitz global optimization
problem (\ref{problem1}) are adopted hereinafter for solving the
problem (\ref{problem2}).  A general scheme describing these methods
is presented in the next section.

\section{A General Scheme describing geometric and information algorithms}
Methods studied in this paper have  a similar structure and belong
to the class of ``Divide-the-best'' global optimization algorithms
introduced in \cite{Divide_the_Best}. They can have the following
differences in their computational schemes distinguishing one
algorithm from another:
\begin{description}
\item(i) Methods are either Geometric or Information (see \cite{Sergeyev&Kvasov(2017),Strongin(1978),
Strongin&Sergeyev(2000)} for   detailed descriptions of these
classes of methods);
\item(ii) Methods can use different approaches for estimating the
 Lipschitz constant: an a priori estimate, a global adaptive estimate, and two
local tuning techniques: Maximum Local Tuning (MLT) and
Maximum-Additive Local Tuning (MALT) (see
\cite{Sergeyev&Kvasov(2017),Sergeyev:et:al.(2016a),Strongin&Sergeyev(2000)}
for   detailed descriptions of these approaches).
\end{description}

The first difference, (i), consists of the choice of characteristics
$R_i$ for the subintervals $[x_{i-1},x_i],~2 \leq i \leq k$, where the
points $x_i,~1 \leq i \leq k$, are called \textit{trial points} and are
points where the objective function $g(x)$ has been evaluated during
previous iterations:
\begin{equation}
\hspace{-2mm}R_i = \left\{
\begin{array}{lr}
\frac{z_i + z_{i-1}}{2} - l_i \frac{x_i - x_{i-1}}{2},&\mbox{\small for
geometric methods},\\
2(z_i + z_{i-1}) - l_i(x_i - x_{i-1}) - \frac{(z_i -
z_{i-1})^2}{l_i (x_i - x_{i-1})},&\mbox{\small for information
methods},
\end{array}
\right.
\label{char}
\end{equation}
where $z_i = g(x_i)$ and $l_i$ is an estimate of the Lipschitz
constant for the subinterval $[x_{i-1},x_i],~2 \leq i \leq k$.

The second distinction, (ii), is related to four different
strategies used to estimate the Lipschitz constant $L$. The first
one consists of applying an a priori given estimate $\overline{L} >
L$. The second way is to use an adaptive global estimate of  the
Lipschitz constant $L$ during the search (the word \textit{global}
means that the same estimate is used for the whole region $D$). The
global adaptive estimate $\overline{L}_k$ can be calculated as
follows
\begin{equation}
\overline{L}_k = \left\{
\begin{array}{lr}
r\cdot H^k,&\mbox{if}~H^k>0,\\
1,&\mbox{otherwise},
\end{array}
\right.
\label{GlobalL}
\end{equation}
where $r>0$ is a reliability parameter and
 \begin{equation}
H^k = \max\{H_i: 2 \leq i \leq k\},
\label{GlobalH}
\end{equation}
\begin{equation}
H_i = \frac{|z_i - z_{i-1}|}{x_i - x_{i-1}},~2 \leq i \leq k.
\label{Hi}
\end{equation}

Finally, the Maximum (MLT) and Maximum-Additive (MALT) local
tuning techniques consist of   estimating  local Lipschitz
constants $l_i$ for each subinterval $[x_{i-1},x_i], 2 \leq i \leq k,$
as follows
\begin{equation}
l_i^{MLT} = \left\{
\begin{array}{lr}
r \cdot \max\{\lambda_i,\gamma_i\}, &\mbox{if}~H^k>0,\\
1,&\mbox{otherwise},
\end{array}
\right.
\label{LocalLM}
\end{equation}

\begin{equation}
l_i^{MALT} = \left\{
\begin{array}{lr}
r \cdot \max\{H_i,\frac{\lambda_i+\gamma_i}{2}\}, &\mbox{if}~H^k>0,\\
1,&\mbox{otherwise},
\end{array}
\right.
\label{LocalLMA}
\end{equation}
where $H_i$ is from (\ref{Hi}), and $\lambda_i$ and $\gamma_i$ are
calculated as follows
\begin{equation}
  \lambda_i= \max\{H_{i-1}, H_i, H_{i+1}\}, \ \ 2 \leq i \leq k,
  \label{lambda}
 \end{equation}
\begin{equation}
  \gamma_i = H^k \frac{(x_i-x_{i-1})}{ X^{max}},\label{gamma}
 \end{equation}
with $H^k$ from (\ref{GlobalH}) and
 \begin{equation}
  X^{max} = \max \{x_i-x_{i-1}: \ 2 \leq i \leq k \}.\label{Xmax}
 \end{equation}
When $i=2$ and $i=k$ only $H_2,$ $H_3$, and
$H_{k-1}, H_k$, should be considered, respectively, in (\ref{lambda}).

After these preliminary descriptions we are ready to describe  the
General Scheme (GS) of algorithms studied in this paper.

 \bd
  \item{\bf Step 0.} \emph{Initialization}.
  Execute first two trials at the points $a$ and $b$, i.\,e., $x^1 :=
  a$, $z^1:=g(a)$ and $x^2 := b$, $z^2:=g(b)$. Set the iteration counter $k := 2$.
Suppose that $k \geq 2$ iterations of the algorithm have already been executed.
  The iteration $k + 1$ consists of the following steps.
 \ed
 \bd
  \item{\bf Step 1.} \emph{Reordering}. Reorder the points
  $x^1,\ldots,x^k$ (and the corresponding function values $z^1,\ldots,z^k$)
  of previous trials by subscripts so that
  $$
    a = x_1 < \ldots < x_k = b, \hspace*{5mm}
    z_i=g(x_i), \, 1 \leq i \leq k.
  $$
  \item{\bf Step 2.} \emph{Estimates of the Lipschitz constant}.
  Calculate the current estimates~$l_i$ of the Lipschitz constant
  for each subinterval $[x_{i-1},x_{i}]$, $2 \leq i \leq k$, in one of the following ways.
  \bd
  \item{\bf Step 2.1.} \emph{A priori given estimate}.
  Take an a priori given estimate $\overline{L}$ of the Lipschitz constant
  for the whole interval $[a,b]$, i.\,e., set $l_i := \overline{L}$.
  \item{\bf Step 2.2.} \emph{Global estimate}.
  Set $l_i := \overline{L}_k,$ where $\overline{L}_k$ is from (\ref{GlobalL}).
  \item{\bf Step 2.3.} \emph{``Maximum'' local tuning}.
  Set $l_i := l_i^{MLT},$ where $l_i^{MLT}$ is from~(\ref{LocalLM}).
  \item{\bf Step 2.4.} \emph{``Maximum-Additive'' local tuning}.
  Set $l_i := l_i^{MALT},$ where $l_i^{MALT}$ is from~(\ref{LocalLMA}).
  \ed
  \item{\bf Step 3.} \emph{Calculation of characteristics}.
  Compute for each subinterval $[x_{i-1},x_i]$, $2 \leq i \leq k$, its characteristic $R_i$
  by using one of the following rules.
  \bd
  \item{\bf Step 3.1.} \emph{Geometric methods}.
   \begin{equation}
    R_i = \frac{z_i + z_{i-1}}{2} - l_i \frac{x_i - x_{i-1}}{2}.
    \label{geomR}
   \end{equation}
  \item{\bf Step 3.2.} \emph{Information methods}.
   \begin{equation}
    R_i = 2(z_i + z_{i-1}) - l_i(x_i - x_{i-1}) -
          \frac{(z_i - z_{i-1})^2}{l_i (x_i - x_{i-1})}.
    \label{infR}
   \end{equation}
  \ed
  \item{\bf Step 4.} \emph{Interval selection}.
  Determine an interval $[x_{t-1},x_t]$, $t=t(k)$, for performing the
  next trial as follows
  \begin{equation} t = \min \arg \min_{2 \leq i \leq k} R_i.
  \label{IntSelection}\end{equation}
  \item{\bf Step 5.} \emph{Stopping rule}. \textbf{If}
    \begin{equation}
      x_t - x_{t-1} \leq \varepsilon,
      \label{epsilon}
    \end{equation}
 where $\varepsilon>0$ is a given accuracy of the global search,
 \textbf{then} \textbf{Stop} and take as an estimate of the global minimum $g^*$
 the value $g_k^* = \min_{1 \leq i \leq k} \{z_i\}$ obtained at a point
 $x_k^* = \arg \min_{1 \leq i \leq k} \{z_i\}$.

 \textbf{Otherwise}, go to Step 6.

 \item{\bf Step 6.} \emph{New trial}. Execute the next trial $z^{k+1}:=g(x^{k+1})$ at the
 point \begin{equation}
 x^{k+1} = \frac{x_t + x_{t-1}}{2} - \frac{z_t - z_{t-1}}{2l_t}.
 \label{newTrial}
 \end{equation}
 Increase the iteration counter $k:=k+1$, and go to Step 1.
 \ed

\section{Strong homogeneity of  algorithms belonging to GS
for finite, infinite,\\  and infinitesimal scailing/shifting
constants}

In this section, we study the strong homogeneity of  algorithms
described in the previous section. This study is executed
simultaneously in the traditional and in the Infinity Computing
frameworks. In fact, so far,  whether these methods were strongly
homogeneous or not was an open problem even for finite constants
$\alpha$ and $\beta$. In this section, we show that methods
belonging to GS enjoy the strong homogeneity property for finite,
infinite, and infinitesimal scaling and shifting constants. Recall
that all the values related to the function $g(x)$ are indicated
by~``~$\widehat{ }$~'' and the values related to the function $f(x)$
are written without ``~$\widehat{ }$~''.

The following lemma shows  how the adaptive estimates of the
Lipschitz constant $\widehat{\overline{L}}_k$,
$\widehat{l}_i^{MLT}$, and $\widehat{l}_i^{MALT}$ that can assume
finite, infinite, and infinitesimal values  are related to the
respective original estimates   $\overline{L}_k$, ${l}_i^{MLT}$, and
${l}_i^{MALT}$ that can be finite only.

\begin{lemma}
 Let us consider the function $g(x) = \alpha f(x) + \beta$, where
$f(x)$ assumes only finite values and has a finite Lipschitz
constant $L$ over the interval $[a,b]$ and $\alpha,~\alpha>0,$ and
$\beta$ can be finite, infinite and infinitesimal numbers. Then,
the adaptive estimates $\widehat{\overline{L}}_k$,
$\widehat{l}_i^{MLT}$ and $\widehat{l}_i^{MALT}$ from  (\ref{GlobalL}), (\ref{LocalLM}) and (\ref{LocalLMA}) are equal
to $\alpha \overline{L}_k$, $\alpha l_i^{MLT}$ and $\alpha l_i^{MALT}$, respectively, if $H^k
> 0$, and to $1$, otherwise.
\end{lemma}
\emph{Proof.} It follows from (\ref{Hi}) that
\begin{equation}
\widehat{H}_i = \frac{|\widehat{z}_i - \widehat{z}_{i-1}|}{x_i - x_{i-1}} = \frac{\alpha |z_i - z_{i-1}|}{x_i - x_{i-1}} = \alpha H_i.
\label{Lemma2:Hi}
\end{equation}

If $H^k \neq 0$, then $H^k = \underset{2 \leq i \leq k}{\max}
\frac{|z_i - z_{i-1}|}{x_i - x_{i-1}}$ and  $H^k \geq H_i,~2 \leq i \leq k$. Thus, using (\ref{Lemma2:Hi}) we  obtain $\alpha H^k
\geq \alpha H_i = \widehat{H}_i$, and, therefore, $\widehat{H}^k =
\alpha H^k$ and from (\ref{GlobalL}) it follows
$\widehat{\overline{L}}_k = \alpha \overline{L}_k$. On the other
hand,  if $H^k = 0$, then both estimates for the functions $g(x)$
and $f(x)$ are equal to 1 (see (\ref{GlobalL})).

The same reasoning can be used to show the respective results for
the local tuning techniques MLT and MALT (see (\ref{LocalLM}) and (\ref{LocalLMA}))
$$\widehat{\lambda}_i = \max\{\widehat{H}_{i-1},\widehat{H}_i,\widehat{H}_{i+1}\} = \alpha \max\{H_{i-1},H_i,H_{i+1}\},$$
$$\widehat{\gamma}_i = \widehat{H}^k \frac{x_i - x_{i-1}}{X^{max}} = \alpha H^k \frac{x_i - x_{i-1}}{X^{max}} = \alpha \gamma_i,$$
$$
\widehat{l}_i^{MLT} = \left\{
\begin{array}{lr}
r \cdot \max\{\widehat{\lambda}_i,\widehat{\gamma}_i\},
&\mbox{if}~\widehat{H}^k>0,\\
1,&\mbox{otherwise}.
\end{array}
\right.
$$
$$
\widehat{l}_i^{MALT} = \left\{
\begin{array}{lr}
r \cdot \max\{\widehat{H}_i,\frac{\widehat{\lambda}_i+\widehat{\gamma}_i}{2}\},
&\mbox{if}~\widehat{H}^k>0,\\
1,&\mbox{otherwise}.
\end{array}
\right.
$$

Therefore, we can conclude that
$$\hspace*{26mm}\widehat{l}_i^{\{MLT,MALT\}} = \left\{ \begin{array}{lr}
\alpha l_i^{\{MLT,MALT\}},& \mbox{if}~H^k>0,\\
1,&\mbox{otherwise}.
\end{array}
\right. \hspace*{18mm}\qed
$$

\begin{lemma}
Suppose that   characteristics $\widehat{R}_i,~2 \leq i \leq k,$ for the
scaled objective function $g(x)$ are equal to an affine
transformation of the characteristics $R_i$ calculated for the
original objective function $f(x)$
 \begin{equation}
  \widehat{R}_i = \widehat{\alpha}_k R_i +
\widehat{\beta}_k,\hspace{5mm} 2 \leq i \leq k, \label{scaledChar}
 \end{equation}
  where scales
$\widehat{\alpha}_k,~\widehat{\alpha}_k>0,$ and $\widehat{\beta}_k$
can be finite, infinite, or infinitesimal and  possibly different
for different iterations $k$. Then, the same interval
$[x_{t-1},x_t]$, $t=t(k)$, from (\ref{IntSelection}) is selected at
each iteration for the next subdivision  during optimizing   $f(x)$
and $g(x)$, i.e., it follows $\widehat{t(k)} = t(k)$.
\end{lemma}
\emph{Proof.} Since due to (\ref{IntSelection})   $t = \arg
\min_{2 \leq i \leq k} R_i$, then $R_t \leq R_i$ and
 $$\widehat{\alpha}_k R_t +\widehat{\beta}_k \leq \widehat{\alpha}_k
R_i + \widehat{\beta}_k,\hspace{5mm}2 \leq i \leq k.
 $$
  That, due to (\ref{scaledChar}), can be re-written as
  $$\widehat{R}_t = \min_{2 \leq i \leq k} \widehat{R}_i =
\widehat{\alpha}_k R_t + \widehat{\beta}_k.
 $$
  Notice that if there are
several values $j$ such that $R_j = R_t$, then (see
(\ref{IntSelection})) we have $t < j, j \neq t$, i.e., even in this
situation it follows $\widehat{t(k)} = t(k)$. This observation
concludes the proof. \hfill $\qed$

The following Theorem shows that methods belonging to the GS enjoy
the strong homogeneity property.

\begin{theorem} Algorithms belonging to  the GS and applied for
solving the problem~(\ref{problem2}) are strongly homogeneous
  for
finite, infinite, and infinitesimal scales $\alpha>0$ and~$\beta$.
\end{theorem}

\emph{Proof.} Two algorithms optimizing functions $f(x)$ and
$g(x)$ will generate the same sequences of trials  if the
following  conditions hold:
\begin{description}
\item (i) The same interval $[x_{t-1},x_t]$, $t=t(k)$, from (\ref{IntSelection})
is selected at each iteration for the  next subdivision  during
optimizing  functions $f(x)$ and $g(x)$, i.e., it follows
$\widehat{t(k)} = t(k)$.
\item (ii) The next trial at the selected interval $[x_{t-1},x_t]$ is
 performed
 at the same point  during
optimizing  functions $f(x)$ and $g(x)$, i.e., in
  (\ref{newTrial}) it follows
$\widehat{x}^{k+1} = x^{k+1}$.
\end{description}

In order to prove assertions (i) and (ii), let us consider
computational steps of the GS. For both functions, $f(x)$ and
$g(x)$, Steps~0 and~1 of the GS work with the same interval $[a,b]$,
do not depend on the objective function, and, as a result, do not
influence (i) and (ii). Step~2 is a preparative one, it is
responsible for estimating the Lipschitz constants for all the
intervals $[x_{i-1},x_i],~2 \leq i \leq k$ and was studied in Lemmas~1--2
above. Step~3 calculates characteristics of the intervals and,
therefore, is directly related to the assertion~(i). In order to
prove it, we consider computations of
characteristics~$\widehat{R}_i$ for all possible cases of
calculating estimates~$l_i$ during Step~2 and show that there always
possible to indicate constants~$\widehat{\alpha}_k$ and
$\widehat{\beta}_k$ from Lemma~3.

Lemmas 1 and 2 show that for the a priori given finite Lipschitz
constant $L$ for the function $f(x)$ (see Step~2.1) it follows
$\widehat{L} = \alpha L$. For the adaptive estimates of the
Lipschitz constants   for   intervals $[x_{i-1},x_i],$ $2 \leq i \leq k,$
(see (\ref{GlobalL}), (\ref{LocalLM}), (\ref{LocalLMA}) and Steps 2.2 -- 2.4 of the GS)
we have $\widehat{l}_i = \alpha l_i$, if $H^k>0$, and $\widehat{l}_i
= l_i = 1$, otherwise (remind that the latter corresponds to the
situation $z_i = z_1, 1 \leq i \leq k$). Since Step~3 includes substeps
defining information and geometric methods, then the following four
combinations of methods with Lipschitz constant estimates computed
at one of the substeps of Step~2 can take place:

\begin{description}
\item(a)
 The value $\widehat{l}_i = \alpha l_i$ and the  geometric method is used.
From (\ref{geomR}) we obtain
$$\widehat{R}_i = \frac{\widehat{z}_{i-1}+\widehat{z}_i}{2} - \widehat{l}_i \frac{x_i - x_{i-1}}{2} = \alpha (\frac{z_{i-1}+z_i}{2} - l_i \frac{x_i - x_{i-1}}{2}) + \beta = \alpha R_i + \beta.$$
Thus, in this case we have $\widehat{\alpha}_k = \alpha$ and
$\widehat{\beta}_k = \beta$.
\item(b)
 The value $\widehat{l}_i = \alpha l_i$ and the information method is used.
 From (\ref{infR}) we get
$$\widehat{R}_i = 2(\widehat{z}_i + \widehat{z}_{i-1}) - \widehat{l}_i(x_i - x_{i-1}) - \frac{(\widehat{z}_i - \widehat{z}_{i-1})^2}{\widehat{l}_i (x_i - x_{i-1})} = $$
$$2\alpha(z_i + z_{i-1}) + 4\beta - \alpha l_i(x_i - x_{i-1}) - \frac{\alpha^2(z_i - z_{i-1})^2}{\alpha l_i (x_i - x_{i-1})} = \alpha R_i + 4 \beta.$$
Therefore, in this case it follows $\widehat{\alpha}_k = \alpha$ and
$\widehat{\beta}_k = 4 \beta$.
\item(c)
 The value $\widehat{l}_i = l_i = 1$ and the geometric method is considered.
 Since in this case $z_i = z_1, 1 \leq i \leq k,$  then  for the geometric method      (see (\ref{geomR}))
we have
$$
\widehat{R}_i = \frac{\widehat{z}_{i-1}+\widehat{z}_i}{2} -
\widehat{l}_i \frac{x_i - x_{i-1}}{2} = \widehat{z}_1 -
\frac{x_i - x_{i-1}}{2} =
 $$
 $$
 \alpha z_1 + \beta -
\frac{x_i - x_{i-1}}{2} = R_i +\alpha z_1 - z_1 + \beta.
$$
Thus, in this case we have $\widehat{\alpha}_k = 1$ and
$\widehat{\beta}_k = z_1(\alpha  -  1) + \beta$.
\item(d)
The value $\widehat{l}_i = l_i = 1$ and the information method is
used. Then, the characteristics (see (\ref{infR})) are calculated as
follows
$$
\widehat{R}_i = 2(\widehat{z}_i + \widehat{z}_{i-1})
- \widehat{l}_i(x_i - x_{i-1}) - \frac{(\widehat{z}_i
- \widehat{z}_{i-1})^2}{\widehat{l}_i (x_i - x_{i-1})} =
 $$
$$
 4\widehat{z}_1 - (x_i - x_{i-1}) = 4 \alpha z_1 + 4\beta - (x_i - x_{i-1})  =
  R_i + 4 \alpha z_1 - 4z_1 + 4\beta.
 $$
Therefore, in this case it follows $\widehat{\alpha}_k = 1$ and
$\widehat{\beta}_k = 4(  z_1 (\alpha -  1) +  \beta)$.
\end{description}

Let us show now that assertion (ii)  also holds.  Since for both the
geometric and the information approaches the the same formula
(\ref{newTrial})  for computing $x^{k+1}$ is used, we should
consider only two cases related to the estimates of the Lipschitz
constant:
\begin{description}
\item(a) If $\widehat{l}_t = \alpha l_t$, then it follows
$$
\widehat{x}^{k+1} = \frac{x_t + x_{t-1}}{2} -
 \frac{\widehat{z}_t - \widehat{z}_{t-1}}{2\widehat{l}_t} =
  \frac{x_t + x_{t-1}}{2} - \frac{\alpha(z_t - z_{t-1})}{2\alpha l_t} =
  x^{k+1}.
  $$
\item(b)
If $\widehat{l}_t = l_t = 1$, then $z_i = z_1,~1 \leq i \leq k,$   and we
have
$$
\widehat{x}^{k+1} = \frac{x_t + x_{t-1}}{2} - \frac{\widehat{z}_t -
 \widehat{z}_{t-1}}{2\widehat{l}_t} = \frac{x_t + x_{t-1}}{2} = x^{k+1}.
 $$
\end{description}
This result concludes the proof. \hfill $\qed$

\section{Numerical illustrations}

In order to illustrate  the behavior of methods belonging to the GS
in the Infinity Computer framework, the following three algorithms
being examples of concrete implementations of the GS have been
tested:
\begin{itemize}
\item
\textbf{Geom-AL:} Geometric method with an a priori given
overestimate of the Lipschitz constant. It is constructed by using
Steps 2.1 and 3.1 in the GS.
\item
\textbf{Inf-GL:} Information method with the global  estimate of the
Lipschitz constant. It is formed by using Steps 2.2 and 3.2 in the
GS.
\item \textbf{Geom-LTM:} Geometric method with the ``Maximum''
local tuning. It is built by applying Steps 2.3 and 3.1 in the GS.

\end{itemize}

The algorithm Geom-AL has   one    parameter -- an a priori given
overestimate of the Lipschitz constant.    In algorithms Geom-LTM
and Inf-GL, the Lipschitz constant  is estimated during the search
and   the reliability parameter $r$ is used. In this work, the
values of the Lipschitz constant of the functions $f(x)$ for the
algorithm Geom-AL have been taken from
\cite{Kvasov&Mukhametzhanov(2018)} (and multiplied by $\alpha$ for
the function $g(x)$). The values of the parameter $r$ for the
algorithms Geom-LTM and Inf-GL have been set to $1.1$ and $1.5$,
respectively. The value $\epsilon = 10^{-4} (b-a)$ has been used in
the stopping criterion (\ref{epsilon}).

Recall that (see Section 2) huge  or very small scaling/shifting
constants can provoke   the ill-conditioning of the scaled function
$g(x)$  in the traditional computational framework. In the Infinity
Computing framework, the positional numeral system
(\ref{grossnumber}) allows us to avoid ill-conditioning and to work
safely with infinite and infinitesimal scaling/shifting constants if
the respective grossdigits and grosspowers are not too large or too
small. In order to illustrate this fact the following two pairs of
the values $\alpha$ and $\beta$ have been used in our experiments:
$(\alpha_1, \beta_1) = (\G1^{-1}, \G1)$ and $(\alpha_2,\beta_2) =
(\G1, \G1^2)$. The corresponding  grossdigits and grosspowers
involved in their representation are, respectively: 1 and $-1$ for
$\alpha_1$; 1 and $1$ for $\beta_1$; 1 and $1$ for   $\alpha_2$; and
1 and $2$ for $\beta_2$. It can be seen that all of these constants
are numbers that do not provoke instability in numerical operations.
Hereinafter scaled functions constructed using constants $(\alpha_1,
\beta_1)$ are indicated as $g(x)$ and functions using $(\alpha_2,
\beta_2)$ are designated as $h(x)$.

\begin{figure}[t!]
\includegraphics[width = 1\linewidth,viewport=0 200  600 600]{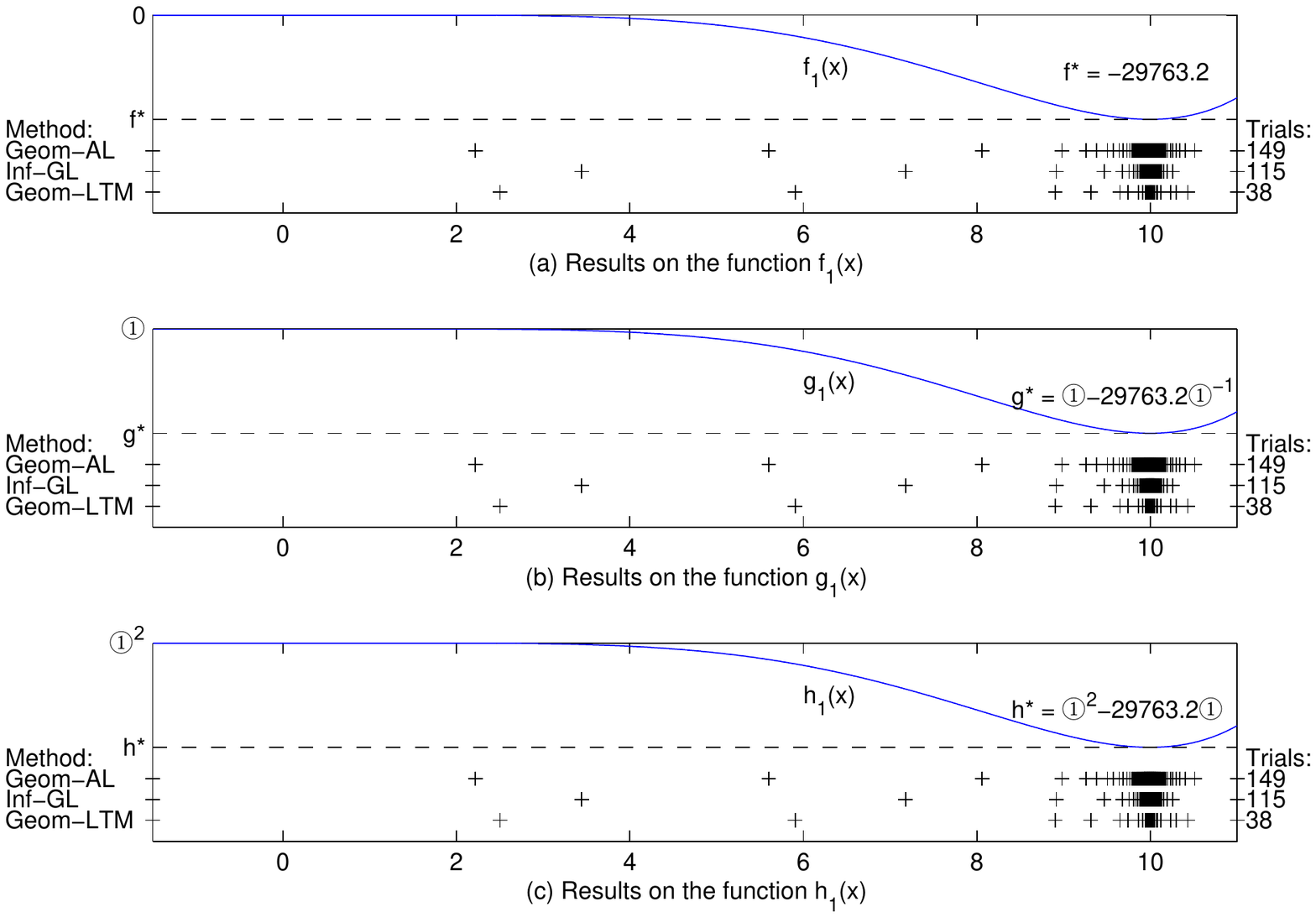}
\caption{Results for (a) the original test function $f_1(x)$ from
\cite{Hansen&Jaumard(1995),Kvasov&Mukhametzhanov(2018)}, (b) the
scaled test function $g_1(x) = \G1^{-1} f_1(x) + \G1$, (c) the
scaled test function $h_1(x) = \G1 f_1(x) + \G1^2$.  Trials are
indicated by the signs ``+'' under the graphs of the functions and
the number of trials for each method is indicated on the right. The
results coincide for each method on all three test functions.}
\label{fig:test0}
\end{figure}
\begin{figure}[t!]
\includegraphics[width = 1\linewidth,viewport=0 200  600 600]{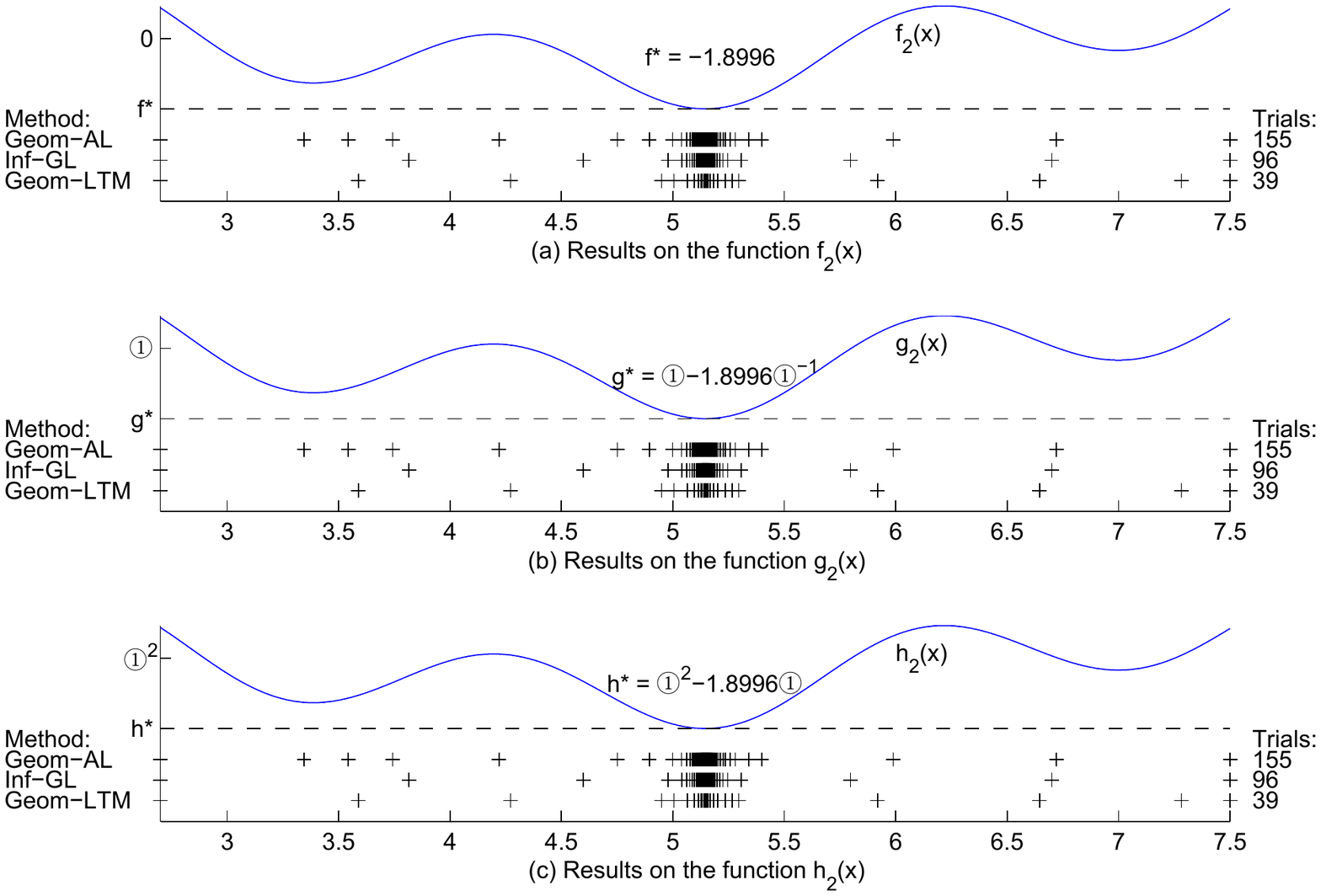}
\caption{Results for (a) the original test function $f_2(x)$ from
\cite{Hansen&Jaumard(1995),Kvasov&Mukhametzhanov(2018)}, (b) the
scaled test function $g_2(x) = \G1^{-1} f_2(x) + \G1$, (c) the
scaled test function $h_2(x) = \G1 f_2(x) + \G1^2$.  Trials are
indicated by the signs ``+'' under the graphs of the functions and
the number of trials for each method is indicated on the right. The
results coincide for each method on all three test functions.}
\label{fig:test1}
\end{figure}

\begin{figure}[t!]
\includegraphics[width = 1\linewidth,viewport=0 200  600 600]{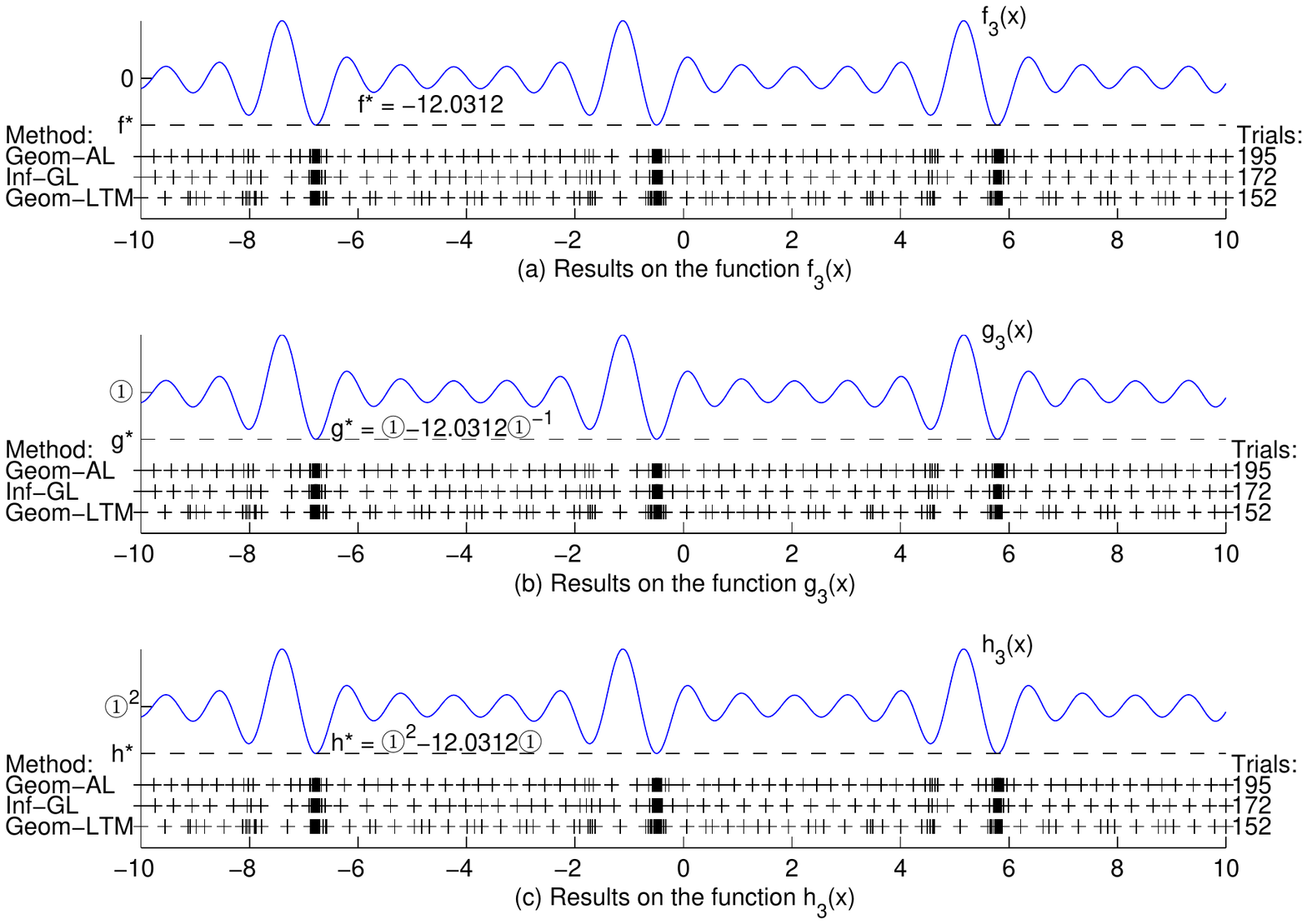}
\caption{Results for (a) the original test function $f_3(x)$ from
\cite{Hansen&Jaumard(1995),Kvasov&Mukhametzhanov(2018)}, (b) the
scaled test function $g_3(x) = \G1^{-1} f_3(x) + \G1$, (c) the
scaled test function $h_3(x) = \G1 f_3(x) + \G1^2$. The results
coincide for each method on all three test functions. The number of
trials for each method is indicated on the right.} \label{fig:test2}
\end{figure}

The algorithms Geom-AL, Inf-GL, and Geom-LTM have been tested on 20
global optimization problems from
\cite{Hansen&Jaumard(1995),Kvasov&Mukhametzhanov(2018)} and on the
respective  scaled functions $g(x)$ and   $h(x)$ constructed from
them. It has been obtained that on all 20 test problems with
  infinite and infinitesimal  constants $(\alpha_1, \beta_1)$
  and   $(\alpha_2, \beta_2)$
the results on the original functions $f(x)$ from
\cite{Hansen&Jaumard(1995),Kvasov&Mukhametzhanov(2018)} and on
scaled functions $g(x)$ and   $h(x)$ coincide. To illustrate this
fact, let us consider the first three problems from the set of 20
tests (see Fig.~\ref{fig:test0}.a, Fig.~\ref{fig:test1}.a, and
Fig.~\ref{fig:test2}.a). They  are defined as follows
 $$f_{1}(x) = \frac{1}{6}x^6 - \frac{52}{25}x^5 +
 \frac{39}{80}x^4+\frac{71}{10}x^3 - \frac{79}{20}x^2 - x + \frac{1}{10},$$
 $$f_{2}(x) =
sin(x)+sin\frac{10x}{3},$$
 $$f_3(x) = \sum_{k=1}^5 -k\cdot
sin[(k+1)x+k].$$
  In Fig.~\ref{fig:test0}.b, Fig.~\ref{fig:test1}.b, and
Fig.~\ref{fig:test2}.b, the results for the scaled functions
$$
g_{i}(x) = \G1^{-1} f_{i}(x) + \G1, \hspace{5mm} i=1,2,3,
$$
 are presented and in Fig.~\ref{fig:test0}.c, Fig.~\ref{fig:test1}.c, and
Fig.~\ref{fig:test2}.c, the results for the scaled functions
 $$h_{i}(x) = \G1 f_{i}(x) + \G1^2, \hspace{5mm} i=1,2,3,
$$
 are shown.
It can be seen that the results coincide for all three methods on
all three test functions $ f_{i}(x), g_{i}(x),$ and $ h_{i}(x),
i=1,2,3$. Analogous results hold for the remaining test problems
from \cite{Hansen&Jaumard(1995),Kvasov&Mukhametzhanov(2018)}.

 In particular, it can be seen from these experiments
that even if the scaling constants $\alpha$ and $\beta$ have a
different order (e.g., when $\alpha$ is infinitesimal and $\beta$ is
infinite) the scaled problems continue to be well-conditioned (cf.
discussion on ill-conditioning in the traditional   framework with
finite scaling/shifting constants, see Fig.~\ref{fig:f}). This fact
suggests that even  if
  finite constants of significantly different orders are required,
$\G1$ can also be used to avoid the ill-conditioning by substituting
very small constants by $ \G1^{-1}$ and very huge constants by
$\G1$. In this case, if, for instance, $\alpha$ is too small (as,
e.g., in (\ref{g3}), $\alpha = 10^{-17}$) and $\beta$ is too large
(as, e.g., in (\ref{g3}), $\beta = 1 \gg 10^{-17}$), the values
$\alpha_1 = \G1^{-1}$ and $\beta_1 = \G1$ can be used in
computations instead of $\alpha =10^{-17}$ and $\beta =1$ avoiding
so underflows and overflows. After the conclusion of the
optimization process, the global minimum of the original function
$f^*$ can be easily extracted from the   solution $g^* = \alpha_1
f^* + \beta_1 = \G1^{-1} f^* +\G1$ of the scaled problem using
$\G1^{-1}$ and $\G1$ and the original finite constants $\alpha$ and
$\beta$ can be used to get the required value $g^* = \alpha f^* +
\beta$ (in our case, $g^* = 10^{-17} f^* + 1$).

\section{Concluding remarks}
Univariate Lipschitz global optimization  problems have been
considered in this paper. Strong homogeneity of   global
optimization algorithms has been studied in the new computational
framework -- Infinity Computing. A new class of global optimization
problems has been introduced  where the objective function can have
finite, infinite or infinitesimal Lipschitz constants. The strong
homogeneity of a class of  geometric and information algorithms used
for solving the univariate Lipschitz global optimization problems
belonging to the new class has been proved for finite, infinite, and
infinitesimal scaling constants. Numerical experiments executed on a
set of test problems taken from the literature confirm the obtained
theoretical results.

Moreover, it has been shown that in cases where global optimization
problems become ill-condi\-tion\-ed in the traditional computational
framework working with finite numbers due to very huge and/or small
scaling/shifting constants, applying Infinity Computing can help in
certain cases. In this situation it is useful to substitute finite
constants provoking problems by infinite and infinitesimal  numbers
that allow one to avoid ill-conditioning of the scaled problems.

\section*{Acknowledgement}
This research   was supported by the Russian Science Foundation,
 project No 15-11-30022 ``Global optimization, supercomputing computations,
 and applications''. The authors thank the unknown reviewers for
 their very useful suggestions.


\begin{thebibliography}{41}
\expandafter\ifx\csname
natexlab\endcsname\relax\def\natexlab#1{#1}\fi
\providecommand{\url}[1]{\texttt{#1}} \providecommand{\href}[2]{#2}
\providecommand{\path}[1]{#1} \providecommand{\DOIprefix}{doi:}
\providecommand{\ArXivprefix}{arXiv:}
\providecommand{\URLprefix}{URL: }
\providecommand{\Pubmedprefix}{pmid:}
\providecommand{\doi}[1]{\href{http://dx.doi.org/#1}{\path{#1}}}
\providecommand{\Pubmed}[1]{\href{pmid:#1}{\path{#1}}}
\providecommand{\bibinfo}[2]{#2} \ifx\xfnm\relax
\def\xfnm[#1]{\unskip,\space#1}\fi
\bibitem[{Amodio et~al.(2017)Amodio, Iavernaro, Mazzia, Mukhametzhanov and
  Sergeyev}]{ODE_3}
\bibinfo{author}{P.~Amodio}, \bibinfo{author}{F.~Iavernaro},
  \bibinfo{author}{F.~Mazzia}, \bibinfo{author}{M.S. Mukhametzhanov},
  \bibinfo{author}{Y.D. Sergeyev}, \bibinfo{title}{A generalized {T}aylor
  method of order three for the solution of initial value problems in standard
  and infinity floating-point arithmetic}, \bibinfo{journal}{Mathematics and
  Computers in Simulation} \bibinfo{volume}{141} (\bibinfo{year}{2017})
  \bibinfo{pages}{24--39}.
\bibitem[{Caldarola(2018)}]{Caldarola_1}
\bibinfo{author}{F.~Caldarola}, \bibinfo{title}{The {S}ierpinski curve viewed
  by numerical computations with infinities and infinitesimals},
  \bibinfo{journal}{Applied Mathematics and Computation} \bibinfo{volume}{318}
  (\bibinfo{year}{2018}) \bibinfo{pages}{321--328}.
\bibitem[{Cococcioni et~al.(2018)Cococcioni, Pappalardo and {Ya.D.
  Sergeyev}}]{Cococcioni}
\bibinfo{author}{M.~Cococcioni}, \bibinfo{author}{M.~Pappalardo},
  \bibinfo{author}{{Ya.D. Sergeyev}}, \bibinfo{title}{Lexicographic
  multi-objective linear programming using grossone methodology: {T}heory and
  algorithm}, \bibinfo{journal}{Applied Mathematics and Computation}
  \bibinfo{volume}{318} (\bibinfo{year}{2018}) \bibinfo{pages}{298--311}.
\bibitem[{{D'Alotto}(2012)}]{DAlotto}
\bibinfo{author}{L.~{D'Alotto}}, \bibinfo{title}{Cellular automata using
  infinite computations}, \bibinfo{journal}{Applied Mathematics and
  Computation} \bibinfo{volume}{218(16)} (\bibinfo{year}{2012})
  \bibinfo{pages}{8077--8082}.
\bibitem[{{De Cosmis} and {De Leone}(2012)}]{DeLeone}
\bibinfo{author}{S.~{De Cosmis}}, \bibinfo{author}{R.~{De Leone}},
  \bibinfo{title}{The use of grossone in mathematical programming and
  operations research}, \bibinfo{journal}{Applied Mathematics and Computation}
  \bibinfo{volume}{218(16)} (\bibinfo{year}{2012}) \bibinfo{pages}{8029--8038}.
\bibitem[{{De Leone}(2018)}]{DeLeone_2}
\bibinfo{author}{R.~{De Leone}}, \bibinfo{title}{Nonlinear programming and
  grossone: {Q}uadratic programming and the role of constraint qualifications},
  \bibinfo{journal}{Applied Mathematics and Computation} \bibinfo{volume}{318}
  (\bibinfo{year}{2018}) \bibinfo{pages}{290--297}.
\bibitem[{Elsakov and Shiryaev(2010)}]{Elsakov&Shiryaev(2006)}
\bibinfo{author}{S.M. Elsakov}, \bibinfo{author}{V.I. Shiryaev},
  \bibinfo{title}{Homogeneous algorithms for multiextremal optimization},
  \bibinfo{journal}{Computational Mathematics and Mathematical Physics}
  \bibinfo{volume}{50} (\bibinfo{year}{2010}) \bibinfo{pages}{1642--1654}.
\bibitem[{Floudas and Pardalos(2009)}]{Floudas&Pardalos(2009)}
\bibinfo{editor}{C.A. Floudas}, \bibinfo{editor}{P.M. Pardalos} (Eds.),
  \bibinfo{title}{Encyclopedia of Optimization (6 Volumes)},
  \bibinfo{edition}{2nd} ed., \bibinfo{publisher}{Springer},
  \bibinfo{year}{2009}.
\bibitem[{Gaudioso et~al.(2018)Gaudioso, Giallombardo and
  Mukhametzhanov}]{Gaudioso&Giallombardo&Mukhametzhanov(2018)}
\bibinfo{author}{M.~Gaudioso}, \bibinfo{author}{G.~Giallombardo},
  \bibinfo{author}{M.S. Mukhametzhanov}, \bibinfo{title}{Numerical
  infinitesimals in a variable metric method for convex nonsmooth
  optimization}, \bibinfo{journal}{Applied Mathematics and Computation}
  \bibinfo{volume}{318} (\bibinfo{year}{2018}) \bibinfo{pages}{312--320}.
\bibitem[{Gergel et~al.(2016)Gergel, Grishagin and
  Gergel}]{Gergel&Grishagin&Gergel(2016)}
\bibinfo{author}{V.~Gergel}, \bibinfo{author}{V.A. Grishagin},
  \bibinfo{author}{A.~Gergel}, \bibinfo{title}{Adaptive nested optimization
  scheme for multidimensional global search}, \bibinfo{journal}{Journal of
  Global Optimization} \bibinfo{volume}{66} (\bibinfo{year}{2016})
  \bibinfo{pages}{35--51}.
\bibitem[{Grishagin et~al.(2018)Grishagin, Israfilov and
  Sergeyev}]{Grishagin_Israfilov_Sergeyev_2018}
\bibinfo{author}{V.A. Grishagin}, \bibinfo{author}{R.A. Israfilov},
  \bibinfo{author}{Y.D. Sergeyev}, \bibinfo{title}{Convergence conditions and
  numerical comparison of global optimization methods based on dimensionality
  reduction schemes}, \bibinfo{journal}{Applied Mathematics and Computation}
  \bibinfo{volume}{318} (\bibinfo{year}{2018}) \bibinfo{pages}{270--280}.
\bibitem[{Hansen and Jaumard(1995)}]{Hansen&Jaumard(1995)}
\bibinfo{author}{P.~Hansen}, \bibinfo{author}{B.~Jaumard},
  \bibinfo{title}{Lipschitz optimization}, in: \bibinfo{editor}{R.~Horst},
  \bibinfo{editor}{P.M. Pardalos} (Eds.), \bibinfo{booktitle}{Handbook of
  Global Optimization}, volume~\bibinfo{volume}{1}, \bibinfo{publisher}{Kluwer
  Academic Publishers}, \bibinfo{address}{Dordrecht}, \bibinfo{year}{1995}, pp.
  \bibinfo{pages}{407--493}.
\bibitem[{Iudin et~al.(2015)Iudin, Sergeyev and Hayakawa}]{Iudin_2}
\bibinfo{author}{D.I. Iudin}, \bibinfo{author}{Y.D. Sergeyev},
  \bibinfo{author}{M.~Hayakawa}, \bibinfo{title}{Infinity computations in
  cellular automaton forest-fire model}, \bibinfo{journal}{Communications in
  Nonlinear Science and Numerical Simulation} \bibinfo{volume}{20(3)}
  (\bibinfo{year}{2015}) \bibinfo{pages}{861--870}.
\bibitem[{Jones et~al.(1993)Jones, Perttunen and
  Stuckman}]{Jones&Perttunen&Stuckman(1993)}
\bibinfo{author}{D.R. Jones}, \bibinfo{author}{C.D. Perttunen},
  \bibinfo{author}{B.E. Stuckman}, \bibinfo{title}{Lipschitzian optimization
  without the {L}ipschitz constant}, \bibinfo{journal}{Journal of Optimization
  Theory and Applications} \bibinfo{volume}{79} (\bibinfo{year}{1993})
  \bibinfo{pages}{157--181}.
\bibitem[{Kvasov and Mukhametzhanov(2018)}]{Kvasov&Mukhametzhanov(2018)}
\bibinfo{author}{D.E. Kvasov}, \bibinfo{author}{M.S. Mukhametzhanov},
  \bibinfo{title}{Metaheuristic vs. deterministic global optimization
  algorithms: {T}he univariate case}, \bibinfo{journal}{Applied Mathematics and
  Computation} \bibinfo{volume}{318} (\bibinfo{year}{2018})
  \bibinfo{pages}{245--259}.
\bibitem[{Margenstern(2015)}]{Margenstern_3}
\bibinfo{author}{M.~Margenstern}, \bibinfo{title}{Fibonacci words, hyperbolic
  tilings and grossone}, \bibinfo{journal}{Communications in Nonlinear Science
  and Numerical Simulation} \bibinfo{volume}{21(1--3)} (\bibinfo{year}{2015})
  \bibinfo{pages}{3--11}.
\bibitem[{Mockus(1988)}]{Mockus(1988)}
\bibinfo{author}{J.~Mockus}, \bibinfo{title}{Bayesian Approach to Global
  Optimization}, \bibinfo{publisher}{Kluwer Academic Publishers},
  \bibinfo{address}{Dodrecht}, \bibinfo{year}{1988}.
\bibitem[{Paulavi\v{c}ius and {\v
  Z}ilinskas(2014)}]{Paulavicius&Zilinskas(2014)}
\bibinfo{author}{R.~Paulavi\v{c}ius}, \bibinfo{author}{J.~{\v Z}ilinskas},
  \bibinfo{title}{Simplicial Global Optimization}, SpringerBriefs in
  Optimization, \bibinfo{publisher}{Springer}, \bibinfo{address}{New York},
  \bibinfo{year}{2014}.
\bibitem[{{Pint\'{e}r}(2002)}]{Pinter(2002)}
\bibinfo{author}{J.D. {Pint\'{e}r}}, \bibinfo{title}{Global optimization:
  software, test problems, and applications}, in: \bibinfo{editor}{P.M.
  Pardalos}, \bibinfo{editor}{H.E. Romeijn} (Eds.),
  \bibinfo{booktitle}{Handbook of Global Optimization},
  volume~\bibinfo{volume}{2}, \bibinfo{publisher}{Kluwer Academic Publishers},
  \bibinfo{address}{Dordrecht}, \bibinfo{year}{2002}, pp.
  \bibinfo{pages}{515--569}.
\bibitem[{Piyavskij(1972)}]{Piyavskij(1972)}
\bibinfo{author}{S.A. Piyavskij}, \bibinfo{title}{An algorithm for finding the
  absolute extremum of a function}, \bibinfo{journal}{USSR Comput. Math. Math.
  Phys.} \bibinfo{volume}{12} (\bibinfo{year}{1972}) \bibinfo{pages}{57--67}.
\bibitem[{Robinson(1996)}]{Robinson}
\bibinfo{author}{A.~Robinson}, \bibinfo{title}{Non-standard Analysis},
  \bibinfo{publisher}{Princeton Univ. Press}, \bibinfo{address}{Princeton},
  \bibinfo{year}{1996}.
\bibitem[{Sergeyev and Garro(2010)}]{Sergeyev_Garro}
\bibinfo{author}{Y.~Sergeyev}, \bibinfo{author}{A.~Garro},
  \bibinfo{title}{Observability of {T}uring machines: A refinement of the
  theory of computation}, \bibinfo{journal}{Informatica}
  \bibinfo{volume}{21(3)} (\bibinfo{year}{2010}) \bibinfo{pages}{425--454}.
\bibitem[{Sergeyev(1998)}]{Divide_the_Best}
\bibinfo{author}{Y.D. Sergeyev}, \bibinfo{title}{On convergence of ``divide the
  best'' global optimization algorithms}, \bibinfo{journal}{Optimization}
  \bibinfo{volume}{44} (\bibinfo{year}{1998}) \bibinfo{pages}{303 -- 325}.
\bibitem[{Sergeyev(2008)}]{Sergeyev(2008)}
\bibinfo{author}{Y.D. Sergeyev}, \bibinfo{title}{A new applied approach for
  executing computations with infinite and infinitesimal quantities},
  \bibinfo{journal}{Informatica} \bibinfo{volume}{19(4)} (\bibinfo{year}{2008})
  \bibinfo{pages}{567--596}.
\bibitem[{Sergeyev(2009)}]{Dif_Calculus}
\bibinfo{author}{Y.D. Sergeyev}, \bibinfo{title}{Numerical point of view on
  {C}alculus for functions assuming finite, infinite, and infinitesimal values
  over finite, infinite, and infinitesimal domains},
  \bibinfo{journal}{Nonlinear Analysis Series A: Theory, Methods $\&$
  Applications} \bibinfo{volume}{71(12)} (\bibinfo{year}{2009})
  \bibinfo{pages}{e1688--e1707}.
\bibitem[{Sergeyev(2011{\natexlab{a}})}]{Num_dif}
\bibinfo{author}{Y.D. Sergeyev}, \bibinfo{title}{Higher order numerical
  differentiation on the {I}nfinity {C}omputer}, \bibinfo{journal}{Optimization
  Letters} \bibinfo{volume}{5(4)} (\bibinfo{year}{2011}{\natexlab{a}})
  \bibinfo{pages}{575--585}.
\bibitem[{Sergeyev(2011{\natexlab{b}})}]{Biology}
\bibinfo{author}{Y.D. Sergeyev}, \bibinfo{title}{Using blinking fractals for
  mathematical modelling of processes of growth in biological systems},
  \bibinfo{journal}{Informatica} \bibinfo{volume}{22(4)}
  (\bibinfo{year}{2011}{\natexlab{b}}) \bibinfo{pages}{559--576}.
\bibitem[{Sergeyev(2017)}]{Sergeyev(2017)}
\bibinfo{author}{Y.D. Sergeyev}, \bibinfo{title}{Numerical infinities and
  infinitesimals: {M}ethodology, applications, and repercussions on two
  {H}ilbert problems}, \bibinfo{journal}{EMS Surveys in Mathematical Sciences}
  (\bibinfo{year}{2017}). \bibinfo{note}{In press}.
\bibitem[{Sergeyev and Grishagin(1994)}]{Sergeyev_Grishagin_1994}
\bibinfo{author}{Y.D. Sergeyev}, \bibinfo{author}{V.A. Grishagin},
  \bibinfo{title}{A parallel method for finding the global minimum of
  univariate functions}, \bibinfo{journal}{Journal of Optimization Theory and
  Applications} \bibinfo{volume}{80} (\bibinfo{year}{1994})
  \bibinfo{pages}{513--536}.
\bibitem[{Sergeyev and Kvasov(2017)}]{Sergeyev&Kvasov(2017)}
\bibinfo{author}{Y.D. Sergeyev}, \bibinfo{author}{D.E. Kvasov},
  \bibinfo{title}{Deterministic global optimization: An introduction to the
  diagonal approach}, \bibinfo{publisher}{Springer}, \bibinfo{address}{New
  York}, \bibinfo{year}{2017}.
\bibitem[{Sergeyev et~al.(2017)Sergeyev, Kvasov and
  Mukhametzhanov}]{Sergeyev:et:al.(2017b)}
\bibinfo{author}{Y.D. Sergeyev}, \bibinfo{author}{D.E. Kvasov},
  \bibinfo{author}{M.S. Mukhametzhanov}, \bibinfo{title}{Operational zones for
  comparing metaheuristic and deterministic one-dimensional global optimization
  algorithms}, \bibinfo{journal}{Mathematics and Computers in Simulation}
  \bibinfo{volume}{141} (\bibinfo{year}{2017}) \bibinfo{pages}{96 -- 109}.
\bibitem[{Sergeyev et~al.(2016{\natexlab{a}})Sergeyev, Mukhametzhanov, Kvasov
  and Lera}]{Sergeyev:et:al.(2016a)}
\bibinfo{author}{Y.D. Sergeyev}, \bibinfo{author}{M.S. Mukhametzhanov},
  \bibinfo{author}{D.E. Kvasov}, \bibinfo{author}{D.~Lera},
  \bibinfo{title}{Derivative-free local tuning and local improvement techniques
  embedded in the univariate global optimization}, \bibinfo{journal}{Journal of
  Optimization Theory and Applications} \bibinfo{volume}{171}
  (\bibinfo{year}{2016}{\natexlab{a}}) \bibinfo{pages}{186--208}.
\bibitem[{Sergeyev et~al.(2016{\natexlab{b}})Sergeyev, Mukhametzhanov, Mazzia,
  Iavernaro and Amodio}]{ODE_2}
\bibinfo{author}{Y.D. Sergeyev}, \bibinfo{author}{M.S. Mukhametzhanov},
  \bibinfo{author}{F.~Mazzia}, \bibinfo{author}{F.~Iavernaro},
  \bibinfo{author}{P.~Amodio}, \bibinfo{title}{Numerical methods for solving
  initial value problems on the {I}nfinity {C}omputer},
  \bibinfo{journal}{International Journal of Unconventional Computing}
  \bibinfo{volume}{12(1)} (\bibinfo{year}{2016}{\natexlab{b}})
  \bibinfo{pages}{3--23}.
\bibitem[{Sergeyev et~al.(2013)Sergeyev, Strongin and
  Lera}]{Sergeyev:et:al.(2013)}
\bibinfo{author}{Y.D. Sergeyev}, \bibinfo{author}{R.G. Strongin},
  \bibinfo{author}{D.~Lera}, \bibinfo{title}{Introduction to Global
  Optimization Exploiting Space-Filling Curves}, \bibinfo{publisher}{Springer},
  \bibinfo{address}{New York}, \bibinfo{year}{2013}.
\bibitem[{Strongin(1978)}]{Strongin(1978)}
\bibinfo{author}{R.G. Strongin}, \bibinfo{title}{Numerical Methods in
  Multiextremal Problems: In\-for\-ma\-tion--Statistical Algorithms},
  \bibinfo{publisher}{Nauka}, \bibinfo{address}{Moscow}, \bibinfo{year}{1978}.
  \bibinfo{note}{(In Russian)}.
\bibitem[{Strongin and Sergeyev(2014)}]{Strongin&Sergeyev(2000)}
\bibinfo{author}{R.G. Strongin}, \bibinfo{author}{Y.D. Sergeyev},
  \bibinfo{title}{Global Optimization with Non-Convex Constraints: {Sequential}
  and Parallel Algorithms}, \bibinfo{publisher}{Kluwer Academic Publishers,
  Dordrecht}, \bibinfo{address}{3rd ed., Springer, New York},
  \bibinfo{year}{2014}.
\bibitem[{\v{Z}ilinskas and \v{Z}ilinskas(2013)}]{Zilinskas&Zilinskas(2013)}
\bibinfo{author}{A.~\v{Z}ilinskas}, \bibinfo{author}{J.~\v{Z}ilinskas},
  \bibinfo{title}{A hybrid global optimization algorithm for non-linear least
  squares regression}, \bibinfo{journal}{Journal of Global Optimization}
  \bibinfo{volume}{56} (\bibinfo{year}{2013}) \bibinfo{pages}{265--277}.
\bibitem[{Zhigljavsky(2012)}]{Zhigljavsky}
\bibinfo{author}{A.A. Zhigljavsky}, \bibinfo{title}{Computing sums of
  conditionally convergent and divergent series using the concept of grossone},
  \bibinfo{journal}{Applied Mathematics and Computation}
  \bibinfo{volume}{218(16)} (\bibinfo{year}{2012}) \bibinfo{pages}{8064--8076}.
\bibitem[{Zhigljavsky and {\v{Z}ilinskas}(2008)}]{Zhigljavsky&Zilinskas(2008)}
\bibinfo{author}{A.A. Zhigljavsky}, \bibinfo{author}{A.~{\v{Z}ilinskas}},
  \bibinfo{title}{Stochastic Global Optimization},
  \bibinfo{publisher}{Springer}, \bibinfo{address}{New\,York},
  \bibinfo{year}{2008}.
\bibitem[{{\v{Z}ilinskas}(1985)}]{Zilinskas(1985)}
\bibinfo{author}{A.~{\v{Z}ilinskas}}, \bibinfo{title}{Axiomatic
  characterization of a global optimization algorithm and investigation of its
  search strategies}, \bibinfo{journal}{Operations Research Letters}
  \bibinfo{volume}{4} (\bibinfo{year}{1985}) \bibinfo{pages}{35--39}.
\bibitem[{{\v{Z}ilinskas}(2012)}]{Zilinskas(2012)}
\bibinfo{author}{A.~{\v{Z}ilinskas}}, \bibinfo{title}{On strong homogeneity of
  two global optimization algorithms based on statistical models of multimodal
  objective functions}, \bibinfo{journal}{Applied Mathematics and Computation}
  \bibinfo{volume}{218} (\bibinfo{year}{2012}) \bibinfo{pages}{8131--8136}.

\end{thebibliography}
\end{document}